\documentclass[10pt]{article}
\usepackage{amsfonts, amssymb, amsmath, amscd}
\usepackage[dvips]{epsfig,graphics}
\usepackage{latexsym, longtable, tabularx}
\usepackage{epsfig, psfrag}

\newtheorem{thm}{Theorem}
\newtheorem{prop}{Proposition}
\newtheorem{lemma}{Lemma}
\newtheorem{cor}[thm]{Corollary}

\newtheorem{defn}{Definition}
\newtheorem{rmk}{Remark}

\newcommand{\Z}{\mathbb{Z}}
\newcommand{\G}{\mathbb{G}}
\newcommand{\s}{\sigma}
\renewcommand{\r}{\rho}

\newcommand{\Q}{\mathbb{Q}}

\renewcommand{\c}{\cite}
\newcommand{\pf}{{\em Proof: \quad }}
\newcommand{\done}{\hfill $\blacksquare$}

\renewcommand{\H}{\mathbb{H}}

\newcommand{\E}{\mathcal{E}}
\newcommand{\I}{\mathcal{I}}
\newcommand{\B}{\mathcal{B}}
\newcommand{\NB}{\mathcal{N}}
\renewcommand{\P}{\mathcal{P}}

\renewcommand{\d}{\mathcal{D}}
\newcommand{\T}{\mathcal{T}}

\newcommand{\e}{\varepsilon}
\def \BRT {{Bollob\'as--Riordan--Tutte }}

\begin{document}

\title{Quasi-tree expansion for the\\ \BRT polynomial}
\author{
Abhijit  Champanerkar
\footnote{The authors gratefully acknowledge support by the National Science Foundation, and support for the first two authors by PSC-CUNY.}\\
{\em {\small Department of Mathematics, College of Staten Island, City University of New York}}\\ \\
Ilya Kofman
\footnotemark[1]
\\
{\em {\small Department of Mathematics, College of Staten Island, City University of New York}}\\ \\
 Neal Stoltzfus
\footnotemark[1]
\\
{\em {\small Department of Mathematics, Louisiana State University}}}
\date{}
\maketitle

\begin{abstract}
\noindent
Bollob\'as and Riordan introduced a three-variable polynomial
extending the Tutte polynomial to oriented ribbon graphs, which are
multi-graphs embedded in oriented surfaces, such that complementary
regions (faces) are discs.  A quasi-tree of a ribbon graph is a
spanning subgraph with one face, which is described by an ordered
chord diagram.  By generalizing Tutte's concept of activity to
quasi-trees, we prove a quasi-tree expansion of the \BRT polynomial.
\end{abstract}

\section{Introduction}
An {\em oriented ribbon graph} is a multi-graph (loops and multiple
edges allowed) that is embedded in an oriented surface, such that its
complement in the surface is a union of $2$--cells.  The embedding
determines a cyclic order on the edges at every vertex.  Terms for the
same or closely related objects include: combinatorial maps, fat
graphs, cyclic graphs, graphs with rotation systems, and dessins
d'enfant (see \c{LZ, BR1} and references therein).

The Tutte polynomial is a fundamental and ubiquitous invariant of
graphs. Bollob\'as and Riordan \c{BR1} extended the Tutte polynomial
to an invariant of oriented ribbon graphs in a way that takes into
account the topology of the ribbon graph.  In \c{BR2}, they
generalized it to a four-variable invariant of non-orientable ribbon
graphs.  We only consider the \BRT polynomial for the orientable case,
and henceforth all ribbon graphs will be oriented.

The Tutte polynomial can be defined by a state sum over all subgraphs,
by contraction-deletion operations, and by a spanning tree expansion
(see \c{ModernGraph} for a detailed introduction)\footnote{The rank
  polynomial, formulated independently by H. Whitney \c{HW}, 
equals the Tutte polynomial after rescaling.}.
  Tutte's original
definition in \c{Tutte} was the spanning tree expansion, discussed
below, which relies on the concept of {\em activity} of edges with
respect to a spanning tree.  In \c{BR1, BR2} the \BRT polynomial was
shown to satisfy many essential properties of the Tutte polynomial,
including a spanning tree expansion using Tutte's activities.

For planar graphs, a spanning tree is a spanning subgraph whose
regular neighborhood has one boundary component.  For ribbon graphs,
the analogue of a spanning tree is a {\em quasi-tree}, which is a
spanning subgraph with one face, introduced in \c{DFKLS2}.  Just as
the spanning trees of a graph determine many of its important
properties, topological properties of a ribbon graph are determined by
the set of its quasi-trees.  A natural question is whether the \BRT
polynomial has a quasi-tree expansion analogous to the spanning tree
expansion for the Tutte polynomial.

In Section \ref{sec_qt_activity}, we extend Tutte's concept of
activity (with respect to a spanning tree) to {\em activity with
  respect to a quasi-tree} by expressing the quasi-tree as an ordered
chord diagram.  For a genus zero ribbon graph, spanning trees and
quasi-trees coincide, and the two notions of activity are the same.
However, for ribbon graphs of higher genus, spanning trees are a
proper subset of quasi-trees, and the two definitions of activity are
quite distinct (see Remark \ref{rmk1} and Section \ref{sec_example}).

In Section \ref{main}, we give an expansion of the \BRT polynomial
over quasi-trees.  Each term in the expansion is determined by a
particular quasi-tree as a product of factors with a topological
meaning.  In the genus zero case, we recover Tutte's original spanning
tree expansion.  In general, our expansion is different from the
spanning tree expansion given in \c{BR2}.  For example, in the case of
one-vertex ribbon graphs, the spanning tree expansion is the same as
the expansion over all subgraphs, but the quasi-tree expansion has
fewer terms (see Remark \ref{rmk2}).
In addition, we show that a specialization of the \BRT polynomial
gives the number of quasi-trees of every genus. 

Together, Sections \ref{bt} and \ref{sec_pf} prove the main theorem, Theorem \ref{mainthm}.
In Section \ref{sec_example}, we compute the quasi-tree expansion for an example.

\section{Activities with respect to a quasi-tree}
\label{sec_qt_activity}

A ribbon graph $\G$ can be considered both as a geometric and as a
combinatorial object.  Starting from the combinatorial definition, let
$(\s_0,\,\s_1,\,\s_2)$ be permutations of $\{1,\ldots,2n\}$, such that
$\s_1$ is a fixed-point free involution and $\s_0\,\s_1\,\s_2=1$.  We
define the orbits of $\s_0$ to be the vertex set $V(\G)$, the orbits
of $\s_1$ to be the edge set $E(\G)$, and the orbits of $\s_2$ to be
the face set $F(\G)$.  Let $v(\G)$, $e(\G)$ and $f(\G)$ be the numbers
of vertices, edges and faces of $\G$.  The preceding data determine an
embedding of $\G$ on a closed orientable surface, denoted $S(\G)$, as
a cell complex.  The set $\{1,\ldots,2n\}$ can be identified with the
directed edges (or half-edges) of $\G$.  Thus, $\G$ is connected if
and only if the group generated by $\s_0,\,\s_1,\,\s_2$ acts
transitively on $\{1,\ldots,2n\}$.  The genus of $S(\G)$ is called the
genus of $\G$, $g(\G)$.  If $\G$ has $k(\G)$ components,
$2g(\G)=2k(\G)-v(\G)+e(\G)-f(\G)=k(\G)+n(\G)-f(\G)$, where
$n(\G)=e(\G)-v(\G)+k(\G)$ denotes the nullity of $\G$.  Henceforth, we
assume that $\G$ is a connected ribbon graph.
See Table \ref{fig1} for an example of distinct ribbon graphs with the
same underlying graph.

\begin{table} \label{rgexample}
\begin{center}
\begin{tabular}{|c|c|}
\hline 
 & \\
\begin{tabular}{ccc}
\hspace*{0.08in} & \includegraphics[width=1in]{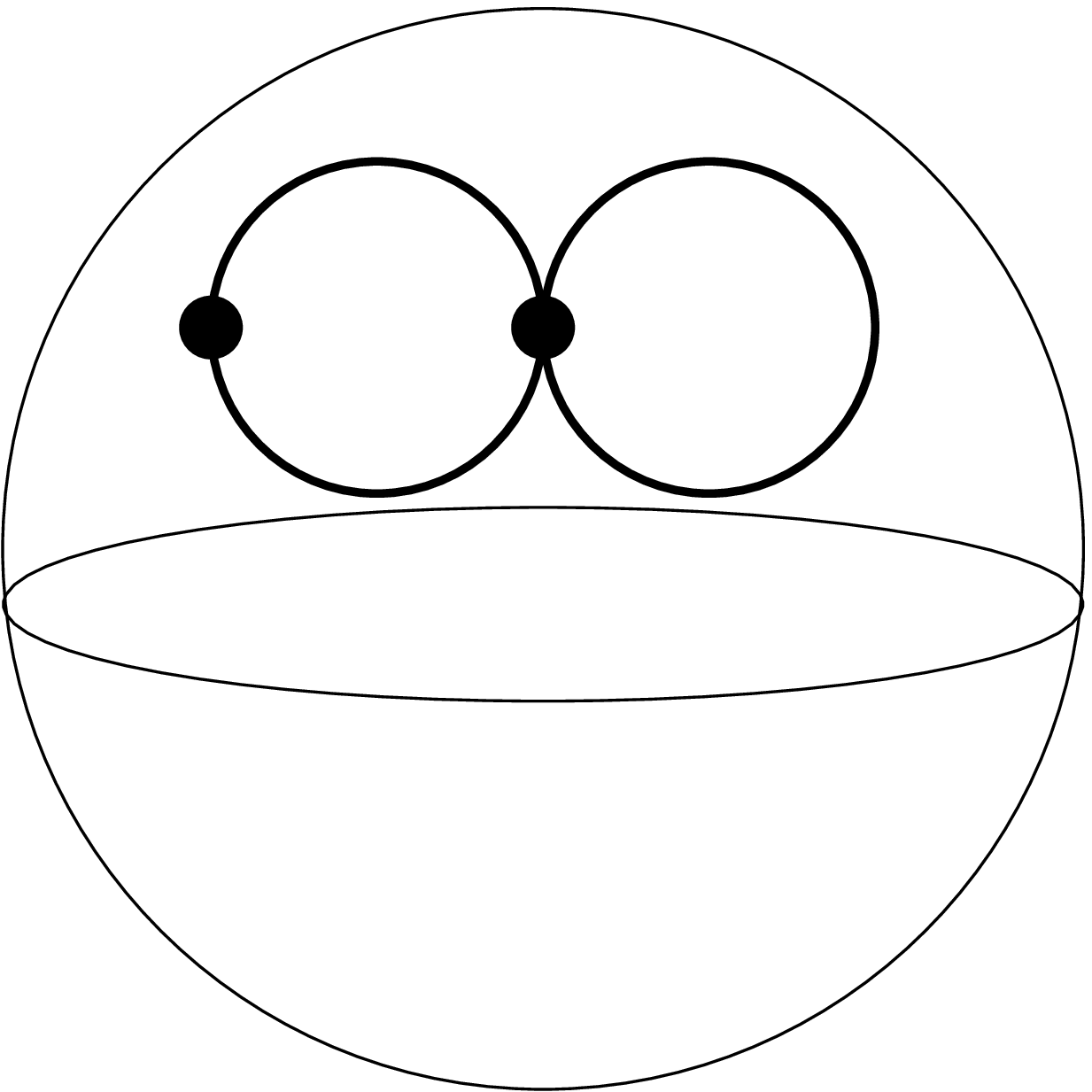} & \hspace*{0.08in} \\
\end{tabular}
&
\begin{tabular}{c}
 \includegraphics[width=1.5in]{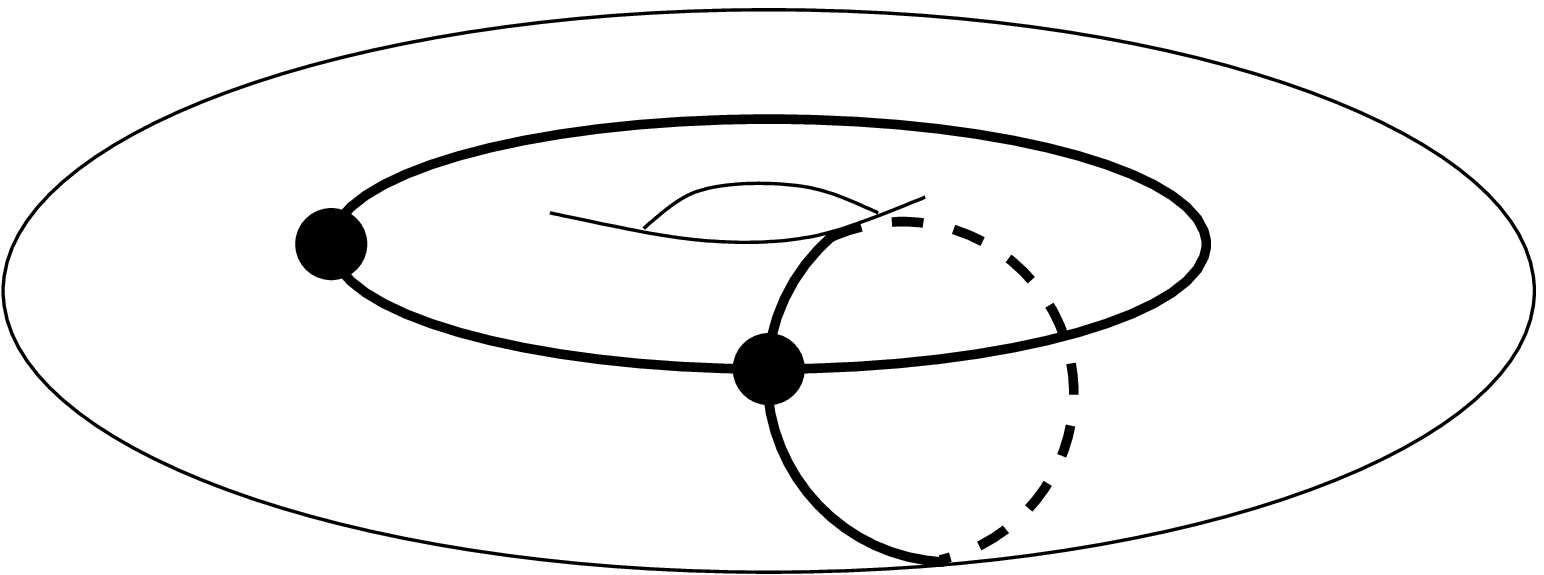} \\
\end{tabular}
\\
\hline
& \\
\psfrag{1}{\scriptsize{$1$}}
\psfrag{2}{\scriptsize{$2$}}
\psfrag{3}{\scriptsize{$3$}}
\psfrag{4}{\scriptsize{$4$}}
\psfrag{5}{\scriptsize{$5$}}
\psfrag{6}{\scriptsize{$6$}}
\includegraphics[width=1.25in]{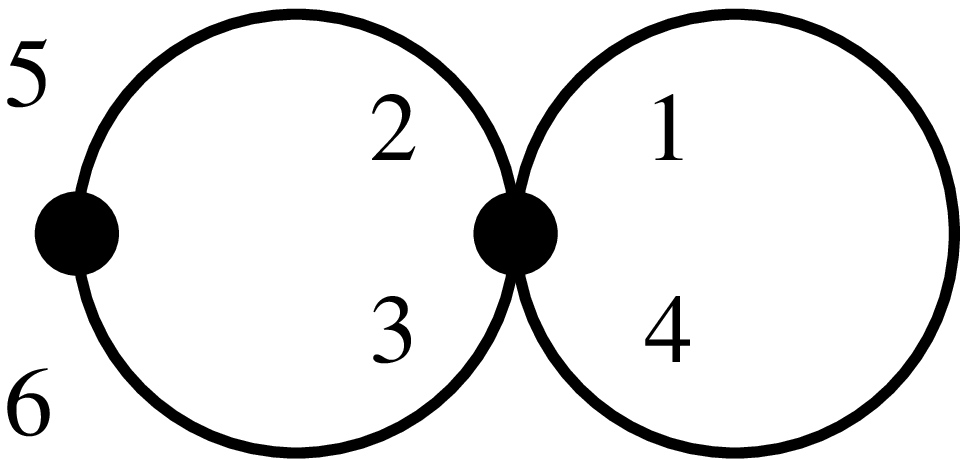} 

& 
\psfrag{1}{\scriptsize{$1$}}
\psfrag{2}{\scriptsize{$2$}}
\psfrag{3}{\scriptsize{$3$}}
\psfrag{4}{\scriptsize{$4$}}
\psfrag{5}{\scriptsize{$5$}}
\psfrag{6}{\scriptsize{$6$}}
\includegraphics[width=1.25in]{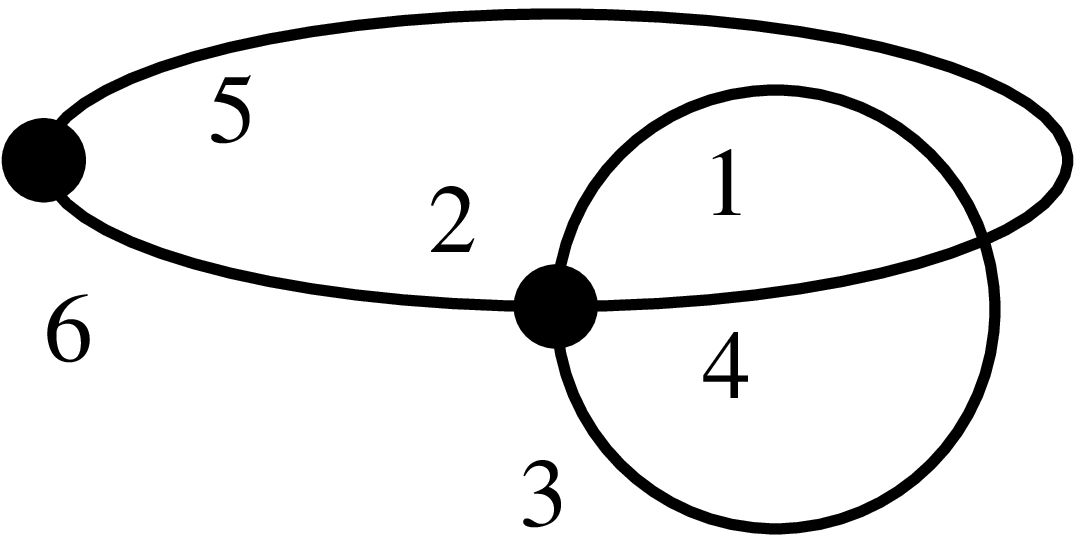}\\
\hline
& \\
\begin{tabular}{l}
 $\sigma_0=(1234)(56)$ \\
 $\sigma_1=(14)(25)(36)$ \\
 $\sigma_2=(246)(35)$ \\
\end{tabular}
&
\begin{tabular}{l}
$\sigma_0=(1234)(56)$ \\
$\sigma_1=(13)(26)(45)$\\
$\sigma_2=(152364)$ \\
\end{tabular}
\\
\hline
\end{tabular}
\end{center}
\caption{Ribbon graphs described as graphs on surfaces and as permutations}
\label{fig1}
\end{table}

Any subgraph $H$ of the underlying graph $G$ of $\G$ determines a
ribbon subgraph $\H$ of $\G$ with underlying graph $H$.  We can
construct its embedding surface $S(\H)$ as follows.
A regular neighborhood of $\H$ can be constructed on the surface
$S(\G)$ by gluing discs at each vertex and rectangular bands whose
midlines are the edges of $\H$.  Let $\gamma_{\H}$ be the union of
simple closed curves that bound such a regular neighborhood of $\H$ on
$S(\G)$.  By attaching a disc to every boundary component of this
regular neighborhood, we construct $S(\H)$, whose genus $g(\H)$
may be smaller than $g(\G)$. By definition, the faces $F(\H)$ are the
complementary regions of $\H$ on $S(\H)$.  Thus, the components of
$\gamma_{\H}$ correspond exactly to the faces $F(\H)$.  So if
$|\gamma_{\H}|$ denotes the number of its components,
$f(\H)=|\gamma_{\H}|$.  In particular, $f(\H)\geq k(\H)$.  Note that
ribbon subgraphs $\H\subseteq\G$ may be disconnected.  Also note that
an isolated vertex $\H$ cannot be represented by
$(\s_0,\,\s_1,\,\s_2)$; in this case, $g(\H)=0$ and $f(\H)=1$.

A ribbon subgraph $\H\subseteq\G$ is called a spanning subgraph if
$V(\H)=V(\G)$.  In this case, $\H$ is a ribbon graph formed from $\G$
by deleting some set of the edges, and keeping all vertices.  
The following concept was introduced and related to the determinant of a link
in \c{DFKLS2}, and also related to Khovanov homology in \c{dkh}.
Following Definition 3.1 of \c{DFKLS2},
\begin{defn}
A \emph{quasi-tree} $\Q$ is a connected spanning subgraph of $\G$ with $f(\Q)=1$.
\end{defn}
Equivalently, a spanning subgraph $\Q$ of $\G$ is a quasi-tree if its
regular neighborhood on $S(\G)$ has exactly one boundary component, $\gamma_{\Q}$.
Also, a spanning connected ribbon graph $\Q$ is a quasi-tree if and
only if $v(\Q)-e(\Q)+2 g(\Q)=1$. If the genus is zero, then the
underlying graph of $\Q$ is a spanning tree.  In Table \ref{fig1},
only the ribbon graph on the right is itself a quasi-tree.

\begin{figure}
\begin{center}
\psfrag{1}{\tiny{$1$}}
\psfrag{2}{\tiny{$2$}}
\psfrag{3}{\tiny{$3$}}
\psfrag{4}{\tiny{$4$}}
\psfrag{5}{\tiny{$5$}}
\psfrag{6}{\tiny{$6$}}
\psfrag{7}{\tiny{$7$}}
\psfrag{8}{\tiny{$8$}}
\psfrag{3c}{\tiny{$3$}}
\psfrag{4c}{\tiny{$4$}}
\psfrag{7c}{\tiny{$7$}}
\psfrag{8c}{\tiny{$8$}}
\includegraphics[height=1.5in]{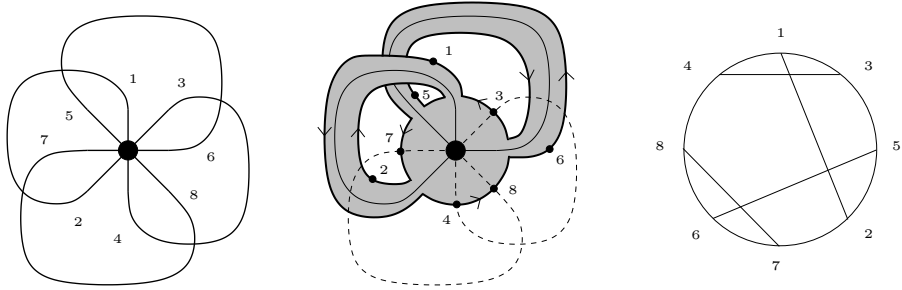}
 \end{center}
 \caption{Ribbon Graph $\G$, quasi-tree $\Q=(12)(56)$ with curve $\gamma_{\Q}$, chord diagram $C_{\Q}$}
 \label{cq}
 \end{figure}

Geometrically, $\gamma_{\Q}$ is a simple closed curve on $S(\G)$ that
divides $S(\G)$ as the connect sum of two surfaces with complementary
genera.  Traversing along $\gamma_{\Q}$, we can mark every half-edge
of $\G$ on its first encounter.  Therefore, $\gamma_{\Q}$ determines
an {\em ordered chord diagram} $C_{\Q}$, which is a circle marked with
$\{1,\ldots,2n\}$ in some order, and chords joining all pairs $\{i,
\s_1(i)\}$.  We say that $\gamma_{\Q}$ is parametrized by $C_{\Q}$.
For example, see Figure \ref{cq}.

\begin{prop}\label{chord}
Let $\G$ be a connected ribbon graph. For every quasi-tree $\Q$ of
$\G$, $\gamma_{\Q}$ is parametrized by the ordered chord diagram
$C_{\Q}$, whose consecutive markings in the positive direction are
given by the permutation:
$$\s(i) =
\begin{cases}
\s_0(i) & i\notin\Q \\
\s_2^{-1}(i) & i\in\Q.
\end{cases} $$
\end{prop}
\pf
Since $\Q$ is a quasi-tree, $\gamma_{\Q}$ is one simple closed curve.
If we choose an orientation on $S(\G)$, we can traverse
$\gamma_{\Q}$ along successive boundaries of bands and vertex discs,
such that we always travel around the boundary of each disc in a
positive direction (i.e., the disc is on the left).  If a half-edge is not in $\Q$, $\gamma_{\Q}$
will pass across it travelling along the boundary of a vertex disc to
the next band. If a half-edge is in $\Q$, $\gamma_{\Q}$ traverses along
one of the edges of its band.
On $\gamma_{\Q}$, we mark a half-edge not in $\Q$ when $\gamma_{\Q}$
passes across it along the boundary of the vertex disc, and we mark a
half-edge in $\Q$ when we traverse an edge of a band in the direction
of the half-edge.  If the half-edge $i$ is not in $\Q$, travelling
along the boundary of a vertex disc, the next half-edge is given by
$\sigma_0$.  If the half-edge $i$ is in $\Q$, traversing the edge of
its band to the vertex disc and then along the boundary of that disc,
the next half-edge is given by $\sigma_0\sigma_1=\sigma_2^{-1}$.  

As $\Q$ is a quasi-tree, each of its half-edges must be in the orbit
of its single face, while the complementary set of half-edges are met
along the boundaries of the vertex discs.
Since we mark all half-edges traversing $\gamma_{\Q}$, the chord diagram $C_{\Q}$ parametrizes $\gamma_{\Q}$.
\done

We now define {\em activity with respect to a quasi-tree}:
\begin{defn}
\label{def_activity}
Fix a total order on the edges of a connected ribbon graph $\G$.
For every quasi-tree $\Q$ of $\G$, this induces an order on the chords
of $C_{\Q}$. A chord is {\em live} if it does not intersect
lower-ordered chords, and otherwise it is {\em dead}.  For any
$\Q$, an edge $e$ is {\em live} or {\em dead} when the
corresponding chord of $C_{\Q}$ is live or dead; and $e$ is {\em
  internal} or {\em external}, according to $e\in\Q$ or $e\in\G - \Q$,
respectively.
\end{defn}

If $\G$ is given by $(\s_0,\,\s_1,\,\s_2)$ as above, we will order the
edges by $\min(i,\s_1(i))$, though any ordering convention will work
as well.  For every quasi-tree $\Q$ of $\G$, the induced order on
chords of $C_{\Q}$ is also given by $\min(i,\s_1(i))$.
In Figure \ref{cq}, we show $C_{\Q}$ such that the only
edge live with respect to $\Q$ is $(12)$, which is internally live.

Tutte \c{Tutte} originally defined activities as follows.  For every
spanning tree $T$ of $G$, each edge $e\in G$ has an activity with
respect to $T$.  If $e \in T$, $\mathit{cut(T,e)}$ is the set of edges
that connect $T\setminus e$.  If $f \notin T$, $\mathit{cyc(T,f)}$ is
the set of edges in the unique cycle of $T \cup f$.  Note $f \in
cut(T,e)$ if and only if $e \in cyc(T,f)$.  An edge $e \in T$ (resp.
$e \notin T$) is {\em internally active} (resp. {\em externally
  active}) if it is the lowest edge in its cut (resp.  cycle), and
otherwise it is {\em inactive}.

Because the two types of activities are distinct, we will use the
notation {\em active/inactive} when referring to activities in the
sense of Tutte with respect to a spanning tree, and {\em live/dead} for activities
with respect to a quasi-tree, as in Definition \ref{def_activity}.

\begin{rmk}\rm
\label{rmk1}
\ 
  \begin{enumerate}
  \item If $g(\G)=0$, then the underlying graph $G$ is planar, and
    $\G$ is given by a fixed planar embedding of $G$.  In this case,
    every quasi-tree $\Q$ of $\G$ is a spanning tree $T$ of $G$.  It
    is easy to check that live (resp. dead) edges of $\G$ with respect
    to $\Q$ are active (resp. inactive) in $G$ with respect to $T$.

  \item A spanning tree of any ribbon graph is also a quasi-tree (of
    genus zero). In this case, the activities using Tutte's original
    definition are different from the activities using our definition.
    For the example in Figure \ref{cq}, the only spanning tree is the
    one with no edges.  Using Tutte's definition, all four edges are
    externally active, but using our definition, the activities are
    $\ell \ell d d$, where $\ell$ and $d$ denote externally live and
    dead, respectively.
    See Section \ref{sec_example} for examples of non-trivial spanning
    trees whose activities are different from those of the corresponding
    quasi-trees.

\item As for planar graphs, the activities with respect to a quasi-tree
   depend on the edge order. In the case of a spanning tree $T$ of a
   planar graph, when the edge order is changed, Tutte proved
   there is a corresponding spanning tree $T'$ whose activity in
the new edge order matches the activity of $T$ in the old order. 
However, for general quasi-trees, such a correspondence may not exist:
In the example in Section \ref{sec_example}, 
switching the edge order by the permutation $(1\,7)\,(2\,8)$
changes the activity of the unique genus 2 quasi-tree from
$LDDDDD$ to $LLLDDD$.

  \end{enumerate}
\end{rmk}

\section{Main results}
\label{main}

The \BRT polynomial $C(\G)\in\Z[X,Y,Z]$ is recursively defined by the disjoint union,
$C(\G_1\amalg\G_2)=C(\G_1)\cdot C(\G_2)$, and the following recursion for edges $e$ of $\G$ and subgraphs $\H$ of $\G$,
where $\G - e$ and $\G/e$ denote deletion and contraction, respectively:
\[ C(\G) =
\begin{cases}
C(\G - e) + C(\G/e) & \text{ if } e \text{ is neither a bridge nor a loop,} \\
X\cdot C(\G/e)       & \text{ if } e \text{ is a bridge,} \\
\sum_{\H} Y^{n(\H)}Z^{g(\H)} & \text{ if } \G \text{ has one vertex,}
\end{cases} \]
where an edge is a {\em bridge} if deleting it increases the number of components. 
Note that $X$ is assigned to a bridge, and $1+Y$ to a loop.  For the
Tutte polynomial $T_G(x,y)$, these are usually $x$ and $y$, respectively.
If $G$ is the underlying graph of a ribbon graph $\G$, then $C(\G;X,Y,1)=T_G(X,1+Y)$.

The \BRT polynomial has a spanning subgraph expansion given by the
following sum over all spanning subgraphs $\H$ of $\G$ (p.85 of
\c{BR2})\footnote{In \c{BR2}, this expansion is given for $R(\G)$.  To
  relate $R(\G)$ to $C(\G)$, we replace $Z$ by $Z^{1/2}$ (p.89 of
  \c{BR2}).},
\begin{equation}
\label{ssdbrt}
 C(\G) = \sum_{\H} (X-1)^{k(\H)-k(\G)}\; Y^{n(\H)}\;Z^{g(\H)}.
\end{equation}

The Tutte polynomial has a spanning tree expansion given by the
following sum over all spanning trees $T$ of a connected graph $G$
with an order on its edges \c{Tutte},
$$ T_G(x,y)= \sum_{T} x^{i(T)}\, y^{j(T)} $$
\noindent where $i(T)$ is the number of internally active edges and
$j(T)$ is the number of externally active edges of $G$ for a given spanning
tree $T$ of $G$.  
Similarly, the \BRT polynomial has the following spanning tree expansion (p.93 of \c{BR2}),
\begin{equation}\label{stbrt}
C(\G)=\sum_{T} X^{i(T)} \sum_{S\subset\e(T)} Y^{n(T\cup S)}\; Z^{g(T\cup S)}
\end{equation}
where $\e(T)$ is the set of externally active edges of $\G$ with respect to a spanning tree $T$ of $\G$.

We will use (\ref{ssdbrt}) to prove a quasi-tree expansion for the
\BRT polynomial, which is different from the expansion (\ref{stbrt}).
Fix a total order on the edges of a connected ribbon graph $\G$.
In Definition \ref{def_activity},
we defined activities ({\em live} or {\em dead}) for edges of $\G$
with respect to $\Q$.  Let $\d(\Q)$ be the spanning subgraph whose edges are the  dead edges in
$\Q$ ({\em internally dead edges}).  Let $\I(\Q)$ be the set of live
edges in $\Q$ ({\em internally live edges}).  Let $\E(\Q)$ be the set
of live edges in $\G-\Q$ ({\em externally live edges}).

For a given quasi-tree let $G_{\Q}$ denote the graph whose vertices
are the components of $\d(\Q)$ and whose edges are the internally live
edges of $\Q$. Let $T_{G_{\Q}}(x,y)$ denote the Tutte polynomial of
$G_{\Q}$.  Our main result is the following:

\begin{thm} \label{mainthm}
Let $\G$ be a connected ribbon graph.  The \BRT polynomial is given by
the following sum over all quasi-trees $\Q$ of $\G$,
$$ C(\G)=\sum_{\Q} Y^{n(\d(\Q))}\;Z^{g(\d(\Q))}\;(1+Y)^{|\E(\Q)|}\;T_{G_{\Q}}(X,1+YZ). $$
\end{thm}

Let $\B(\Q)$ and $\NB(\Q)$ be the set of internally live edges of $\Q$
that are, respectively, bridges and edges that join the same component
of $\d(\Q)$.  Thus, $G_{\Q}$ has $|\B|$ bridges and $|\NB|$ loops,
which contribute factors $X^{|\B|}$ and $(1+YZ)^{|\NB|}$ to
$T_{G_{\Q}}(X,1+YZ)$ in Theorem \ref{mainthm}.

In the case when $\G$ has a single vertex, there are only loops, so we have the following simplification:
\begin{cor} \label{onevertex}
Let $\G$ be a connected ribbon graph with one vertex.
Taking the sum over all quasi-trees $\Q$ of $\G$,
$$ C(\G)=\sum_{\Q} Y^{n(\d(\Q))}\;Z^{g(\d(\Q))}\;(1+Y)^{|\E(\Q)|}\;(1+YZ)^{|\I(\Q)|}.$$
\end{cor}

\begin{rmk}\rm
\label{rmk2}
\
\begin{enumerate}
\item If $g(\G)=0$, by Remark \ref{rmk1}(i), quasi-trees of $\G$ are
  spanning trees of the underlying graph $G$, and live (resp. dead)
  reduces to active (resp.  inactive).  In this case, $G_{\Q}$ is a
  tree with $|\I(\Q)|$ edges.  After substituting $Y=y-1$ and $Z=1$ in
  $C(\G)$, we recover Tutte's original spanning tree expansion for
  $T_G(x,y)$ from Theorem \ref{mainthm}.

\item For one-vertex ribbon graphs, the only spanning tree is the
  subgraph with no edges.  All edges are loops, so all edges are
  externally active in the sense of Tutte. The spanning tree expansion
  (\ref{stbrt}) becomes the expansion (\ref {ssdbrt}) over all
  subgraphs.  In contrast, the quasi-tree expansion in Corollary
  \ref{onevertex} has fewer terms because some subgraphs are not
  quasi-trees.

\item Ed Dewey \cite{Dewey} has generalized both activity with respect to a quasi-tree and Theorem \ref{mainthm} to the non-orientable case.
\footnote{This work was part of the NSF-supported Research Experience for Undergraduates at LSU.}
\end{enumerate}
\end{rmk}

\subsection{Counting quasi-trees}

The Tutte polynomial counts the number of spanning trees of a 
connected graph $G$ by the specialization $T_G(1,1)$. Below, we 
show that specializing the \BRT polynomial counts the 
number of quasi-trees of every genus.

\begin{prop}
\label{qtcount}
Let $q(\G;t,Y)=C(\G;1,Y,t Y^{-2})$. Then $q(\G;t,Y)$ is a polynomial
in $t$ and $Y$ such that 
$$q(\G;t,0)=\sum\nolimits_j a_j t^j$$ 
where $a_j$ is the number of quasi-trees of genus $j$. Consequently, 
 $q(\G;1,0)$ equals the number of quasi-trees of $\G$. 

\end{prop}
\pf The surviving terms in the expansion (\ref{ssdbrt}) of
$C(\G;1,Y,Z)$ satisfy $k(\H)=k(\G)=1$, so they correspond to connected
spanning subgraphs.  Hence,
$$ q(\G;t,Y)=C(\G;1,Y,t Y^{-2})= \sum_{\H} t^{g(\H)}\; Y^{n(\H)-2g(\H)} $$
where the sum is taken over connected spanning subgraphs. 
Since $2g(\H)=k(\H) +n(\H)-f(\H)$, it follows that 
$n(\H)-2g(\H)=f(\H)-k(\H) =f(\H)-1 \geq 0$. This proves 
that $q(\G;t,Y)$ is a polynomial. 
The terms of $q(\G,t,0)$ are those whose $Y$ exponent vanishes, which come from spanning subgraphs $\H$ with $f(\H)=1$.
These are precisely quasi-trees, whose genus $g(\H)$ is given by the exponent on $t$. 
\done

\subsection{Duality}
\label{sec_duality}
The Tutte polynomial satisfies an important duality property; for a
dual graph $G^*$, $T_G(x,y)=T_{G^*}(y,x)$.  Since a ribbon graph is
embedded in a surface, there is a natural dual ribbon graph.
Bollob\'as and Riordan \c{BR2} found a 1--variable specialization of
the \BRT polynomial that is invariant under this duality.

Building on the work of Ellis-Monaghan and Moffat,
Chmutov found the \BRT polynomial satisfies a much more general
duality with respect to any subset of edges of a ribbon graph (see
\c{Chmutov} and references therein).  When all the edges are dualized,
this construction yields the usual dual ribbon graph.  Let $g$ denote
the genus of $\G$.  In our notation, we have
$$ (X-1)^g\, C_{\G}(X,Y,Z)\big |_{(X-1)YZ=1} = Y^g\, C_{\G^*}(Y,X,Z)\big |_{(X-1)YZ=1} $$

More recently, Krushkal \c{Krushkal} introduced a four-variable
polynomial invariant of orientable ribbon graphs that satisfies a
duality relation like the Tutte polynomial, and that specializes to
the \BRT polynomial.

The quasi-trees of a ribbon graph and its dual are in one-one
correspondence.  Since $\gamma_{\Q}$ is a simple closed curve on
$S(\G)$ that divides $S(\G)$ as the connect sum of two surfaces with
complementary genera, the genus of the dual quasi-tree
$g(\Q^*)=g(\G)-g(\Q)$ (Theorem 4.1 of \c{DFKLS2}).  It is an
interesting question to understand the above duality in terms of the
quasi-tree expansion, and whether this expansion gives rise to new
duality properties.

\section{Binary tree of spanning subgraphs} \label{bt}

The spanning subgraphs of a given ribbon graph $\G$ form a poset (of
states) $\P$ isomorphic to the boolean lattice, $\{0,1\}^{E(\G)}$ of
subsets of the set of edges. The partial order is given by $\mathcal
E=(e_i)\preceq \mathcal E'=(e'_i)$ provided $e_i \leq e'_i$ for all
$i$.  In this section, we define a binary tree $\T$, which is similar
to the skein resolution tree for diagrams widely used in knot theory
(see, e.g., \c{KnotsPhysics}). By the construction below, the leaves of $\T$
correspond exactly to quasi-trees of $\G$.

A {\em resolution} of $\G$ is a function $s: E(\G)\to\{0,1\}$, which
determines a spanning subgraph $\H_s=\{e\in\G\ |\ s(e)=1\}$.  Let
$\r:E(\G) \to \{0,1,* \}$ be a {\em partial resolution} of $\G$, with
edges called {\em unresolved} if they are assigned $*$.  Let
$\H_{\r}=\{e\in\G\ |\ \r(e)=1\}$.  A partial resolution determines an
interval in the poset, $[\rho]= \{s\ |\ s(e_i)=\rho(e_i) \
\mathrm{if}\ \rho(e_i)\in\{0,1\}\}=[\rho\wedge 0,\rho\wedge 1]$, the
interval between $\rho\wedge 0$ with all unresolved edges of $\rho$
set to zero, and $\rho\wedge 1$ with all unresolved edges of $\rho$
set to one.
Given a partial resolution $\r$, we call both $\r$ and $\H_{\r}$ {\em split} if
$f(\H_{\r}\cup U) > 1$ for all subsets $U$ of unresolved edges.  

\begin{figure}
\begin{center}
\psfrag{1}{\footnotesize{$1$}}
\psfrag{0}{\footnotesize{$0$}}
\psfrag{s}{\footnotesize{$*$}}
\psfrag{sd}{\footnotesize{split}}
\psfrag{e}{$e$}
\includegraphics[width=5.5in]{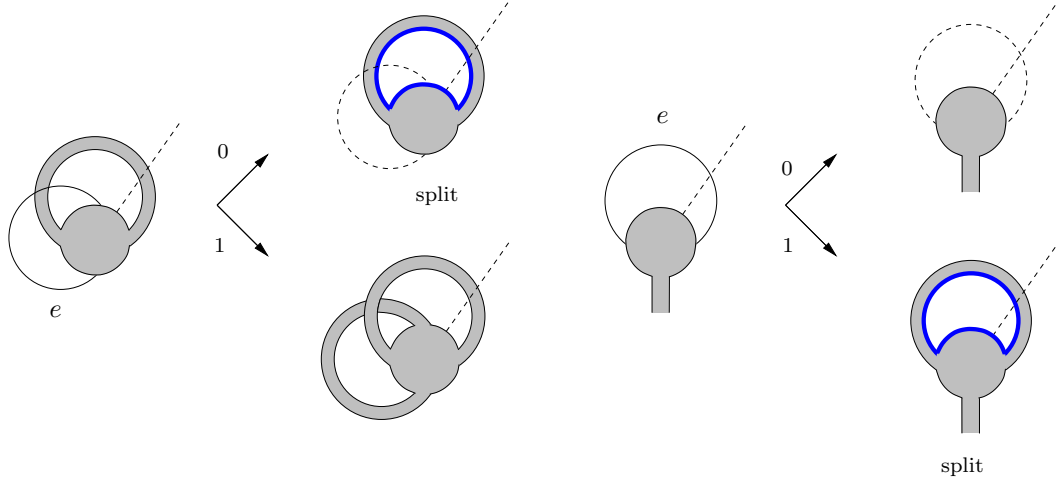}
 \end{center}
 \caption{Two possibilities for a nugatory edge $e$: When $e$ is resolved 
as indicated, the thicker boundary component remains disjoint for 
all choices of unresolved edges, resulting in at least two boundary 
components in the ``split'' cases.}
 \label{nugatoryfig}
\end{figure}

\begin{defn}\label{nugatorydef2}
If $e$ is an unresolved edge in a partial resolution $\rho$, let
$\r^e_0,\, \r^e_1$ be partial resolutions obtained from $\rho$ by
resolving $e$ to be $0$ and $1$, respectively.  Then $e$ is called
{\em nugatory} if either one of $\H_{\r_0}$ or $\H_{\r_1}$ is split.
\end{defn}
Note that an unresolved edge $e$ of $\r$ is nugatory if and only if 
one of the intervals $[\rho^e_0]$ or $[\rho^e_1]$ contains no quasi-trees.
Figure \ref{nugatoryfig} shows two possibilities for a nugatory edge.

 For example, when $g(\G)=0$ and $\r$ is not split, an edge $e$ is
 nugatory in $\r$ if and only if adding it completes a cycle in
 $\rho^e_1$, or $\rho^e_0$ is disconnected and no unresolved edges can
 connect it back.

\begin{thm}\label{treethm}
For any connected ribbon graph $\G$ with ordered edges, there exists a
rooted binary tree $\T$ whose nodes are partial resolutions
$\r$ of $\G$, and whose leaves correspond to quasi-trees $\Q$ of $\G$.
If the leaf $\r$ corresponds to $\Q$, then its unresolved edges are
nugatory, and they can be uniquely resolved to obtain $\Q$.  In $\G$,
these are exactly the live edges with respect to $\Q$.
\end{thm}
\pf 
We prove this theorem in a sequence of two lemmas below.

Let the root of $\T$ be the totally unresolved partial resolution,
$\r(e)=*$ for all $e$.  We resolve edges by changing $*$ to $0$ or $1$
in the reverse order (starting with highest ordered edge). If an edge
is nugatory, the edge is left unresolved, and we proceed to the next
edge.  For a given node $\r$ in $\T$, if $e$ is not nugatory then the
left child is $\r^e_0$ and the right child is $\r^e_1$.  We terminate
this process at a leaf when all subsequent edges are nugatory, and
return as far back up $\T$ as necessary to a node with a non-nugatory
edge still left to be resolved.  Therefore, the leaves of $\T$ are
spanning subgraphs of $\G$ all of whose unresolved edges are nugatory.

Let $\gamma_{\r} = \gamma_{\H_{\r}}$, which was defined previously as
the boundary of a certain regular neighborhood of $\H_{\r}$, and let
$|\gamma_{\r}|$ denote the number of its components.  By definition,
$f(\H_{\r})=|\gamma_{\r}|$, which is the number of faces on $S(\H_{\r})$, the associated
surface for $\H_{\r}$.

Let $\Gamma(\r)=\gamma_{\r} \cup {\rm \ Int}(\r^{-1}(*))$, where ${\rm \ Int}(\r^{-1}(*))$ denotes the set of interiors of all unresolved edges on $S(\G)$.
Note that $\Gamma(\r)$ is connected if and only if we can join the components of $\gamma_{\r}$ by resolving some edges of $\rho$.
Since $f(\H_{\r})=|\gamma_{\r}|$, $\r$ is split if and only if $\Gamma(\r)$ is disconnected.  

\begin{lemma}\label{nugatorylemma}
Let $\r$ be any partial resolution that is not split, with an unresolved edge $e\in\G$.
For $i\in\{0,1\}$, let $\r_i=\rho^e_i$, and let $\Gamma_i(\r,e)=\Gamma(\rho^e_i)$.
The edge $e$ is nugatory if and only if either $\Gamma_0(\r,e)$ or
$\Gamma_1(\r,e)$ is disconnected on $S(\G)$.
If $\Gamma_0(\r,e)$ is disconnected then $|\gamma_{\r_1}| = |\gamma_{\r}|-1$
and $|\gamma_{\r_0}| = |\gamma_{\r}|$.  If $\Gamma_1(\r,e)$ is
disconnected then $|\gamma_{\r_0}| = |\gamma_{\r}|$ and $|\gamma_{\r_1}|
= |\gamma_{\r}| +1$.
\end{lemma}
\pf 
For $i\in\{0,1\}$, $\r_i$ is split if and only if $\Gamma_i(\r,e)$ is disconnected.
Since $\r$ is not split, $\H_{\r_0}$ or $\H_{\r_1}$ is split
if and only if deleting $e$ or cutting along $e$, respectively,
disconnects $\Gamma(\r)$.

If $\Gamma_0(\r,e)$ is disconnected then $e$ is the only edge
connecting two components of $\gamma_{\r}$.  Hence, these two
components are connected in $\gamma_{\r_1}$.
This gives $|\gamma_{\r_1}| = |\gamma_{\r}|-1$ and $|\gamma_{\r_0}| =
|\gamma_{\r}|$.  On the other hand, if $\Gamma_1(\r,e)$ is
disconnected, then $e$ intersects a component of $\gamma_{\r}$ twice
without linking any other unresolved edge, so this component becomes
disconnected in $\gamma_{\r_1}$. This gives $|\gamma_{\r_0}| =
|\gamma_{\r}|$ and $|\gamma_{\r_1}| = |\gamma_{\r}| +1$. 
\done

We can now see that the partial resolution of a leaf can be resolved uniquely to give a quasi-tree.
By construction, for a leaf $\r$ of $\T$, $\H_{\r}$ is not split, so
there exists a resolution $s\in [\r]$ such that
$f(\H_s)=|\gamma_{\H_s}|=1$.  In particular, since all unresolved
edges are nugatory, by Lemma \ref{nugatorylemma}, there is a unique
resolution $s\in [\r]$ such that $|\gamma_{\H_s}|$ is
minimized. Including nugatory edges $e$ for which $\Gamma_1(\r,e)$ is
connected, and excluding nugatory edges $e$ for which $\Gamma_1(\r,e)$
is disconnected, $|\gamma_{\H_s}|=1$.  Hence, $\H_s$ is a quasi-tree.

\begin{lemma} \label{qtactivity} 
  Let $\r$ be a leaf of $\T$, and let $\Q\in [\r]$ be the
  corresponding quasi-tree.  If $\r(e)=*$ then $e$ is live with
  respect to $\Q$, and otherwise $e$ is dead with respect to $\Q$.
\end{lemma}

\pf  If $e_i$ and $e_j$ are any edges of $\r$, we will say that $e_i$
and $e_j$ link each other if, when uniquely resolved to obtain $\Q$, their
endpoints alternate on $\gamma_{\Q}$. Equivalently, their
corresponding chords intersect in $C_{\Q}$.  This notion does not
depend on whether the edges are resolved in $\r$. 
If $g(\G)=0$, then $\Q$ is a spanning tree, and edges link
each other if and only if they satisfy a cut-cycle condition with
respect to $\Q$: $e_i\in cut(\Q,e_j)$ or $e_i\in cyc(\Q,e_j)$.

Let $e_i$ and $e_j$ be unresolved edges of $\r$, which are therefore
nugatory.  Let $s\in [\r]$ be the unique resolution such that
$\H_s=\Q$.  Let $s'$ be the resolution obtained from $s$ by changing
the states of both $e_i$ and $e_j$.  If $e_i$ and $e_j$ link each
other, then $|\gamma_{s'}|=|\gamma_{s}|=1$.
Hence, $\H_{s'}$ is a quasi-tree for a second resolution $s' \in [\r]$,
which is a contradiction.  Thus, unresolved edges can only link
resolved edges.

Suppose $e_i$ is unresolved and links a resolved edge $e_j$ with
$j<i$.  There exists a unique closest parent $\tilde\r$ of $\r$ in
$\T$, such that $e_j$ is a non-nugatory unresolved edge in $\tilde\r$.
Since edges are resolved in the reverse order, $e_i$ is nugatory in
$\tilde\r$.  As $e_i$ links $e_j$, $\Gamma_0(\tilde\r,e_i)$ and
$\Gamma_1(\tilde\r,e_i)$ are both connected, which contradicts Lemma
\ref{nugatorylemma}.  Thus, if $e_i$ and $e_j$ are linked then $i<j$,
so $e_i$ is live.

Now, let $e_i$ be a resolved edge of $\r$.  There exists a unique
closest parent $\tilde\r$ of $\r$ in $\T$, such that $e_i$ is a
non-nugatory unresolved edge in $\tilde\r$.  By Lemma
\ref{nugatorylemma}, $\Gamma_0(\tilde\r,e_i)$ and
$\Gamma_1(\tilde\r,e_i)$ are both connected.  Hence, there exists
$e_j$, which is unresolved in $\tilde\r$, such that $e_i$ and $e_j$
are linked.  If $e_j$ is resolved after $e_i$ in $\T$, then $j<i$.
Since $e_i$ and $e_j$ are linked, $e_i$ is dead.  On the other hand,
if $e_j$ is left unresolved in $\T$, then $e_j$ is live by the
argument in the previous paragraph with $i$ and $j$ reversed.  Since
$e_i$ and $e_j$ are linked, and $e_j$ is live, it follows that $e_i$
is dead.  \done

This completes the proof of Theorem \ref{treethm}. \done

\section{Proof of Theorem \ref{mainthm}}\label{sec_pf}

Let $\H \subseteq \G$ be a spanning subgraph. Let $n(\H)$, $g(\H)$ and $k(\H)$
denote the nullity, genus and number of components of $\H$, respectively. Since $v(\H)=v(\G)$,
$$n(\H)=k(\H) - v(\G)+e(\H),\ \  \ g(\H)=\frac{2k(\H) - v(\G)+e(\H)-f(\H)}{2}.$$
Let $\Q$ be a quasi-tree of $\G$. Let $\I=\I(\Q)$ and $\E=\E(\Q)$ be
the internally and externally live edges with respect to $\Q$. 
Let $\d=\d(\Q)$ be the spanning subgraph whose edges are the  dead edges in
$\Q$.

By Theorem \ref{treethm}, there is a unique partial resolution $\r$ of $\G$ that is a leaf of $\T$,
for which $\Q\in [\r]$, and all resolutions $\H_s$ for $s\in [\r]$ are
of the form $\d\cup S$ where $S \subseteq \I \cup \E$.
All resolutions $\H_s$ are elements of the state poset $\P$, so the sum
in (\ref{ssdbrt}) is a state sum for $\P$.  The sum in Theorem
\ref{mainthm} is a state sum for $\T$.  Below, we prove that these
two state sums are equal.

\begin{lemma} \label{extlivecount} 
For a quasi-tree $\Q$ of $\G$, let $S = S_1 \cup S_2$, where $S_1
\subseteq \I(\Q)$ and $ S_2 \subseteq \E(\Q)$.
\begin{enumerate}
\item $k(\d(\Q)\cup S)=k(\d(\Q)\cup S_1)$
\item $n(\d(\Q)\cup S)=n(\d(\Q)\cup S_1)+ |S_2|$
\item $g(\d(\Q)\cup S)=g(\d(\Q)\cup S_1)$.
\end{enumerate}
\end{lemma}
\pf Let $e\in\E(\Q)$. By Theorem \ref{treethm}, $\Q$ corresponds to
$\r$ such that $e$ is nugatory.  By Lemma \ref{nugatorylemma},
$\Gamma_1(\r,e)$ is disconnected, so
$\Gamma_0(\r,e)$ is connected.  Hence, $e$ intersects only one
component of $\gamma_{\d}$.  Thus, $k(\d\cup e)=k(\d)$, and part 1 follows.
\begin{eqnarray*}
n(\d\cup S) &=& k(\d\cup S) - v(\G)+e(\d\cup S) \\
&=& k(\d\cup S_1)- v(\G) + e(\d\cup S_1)+|S_2| \\
&=& n(\d\cup S_1) +|S_2|.
\end{eqnarray*}
Since $f(\H)=|\gamma_{\H}|$, by Lemma \ref{nugatorylemma}, $f(\d\cup e)=f(\d)+1$, hence
\begin{eqnarray*}
2 g(\d\cup S)&=& 2k(\d\cup S)-v(\G)+ e(\d\cup S) -f(\d\cup S) \\
&=&2k(\d\cup S_1)-v(\G)+\big(e(\d\cup S_1)+|S_2|\big)-\big(f(\d\cup S_1)+|S_2|\big)\\
&=& 2k(\d\cup S_1)-v(\G)+e(\d\cup S_1)-f(\d\cup S_1) \\
&=& 2 g(\d\cup S_1).
\end{eqnarray*}
\done
\begin{lemma} \label{intlivecount} 
For a quasi-tree $\Q$ of $\G$, let $S_1 \subseteq \I(\Q)$.  Let $W$ be
the spanning subgraph of $G_{\Q}$ whose edges are the edges in $S_1$.
\begin{enumerate}
\item $n(\d(\Q)\cup S_1)= n(\d(\Q))+ n(W)$
\item $ g(\d(\Q)\cup S_1)= g(\d(\Q)) + n(W)$.
\end{enumerate}
\end{lemma}
\pf For spanning subgraph $W$ of $G_{\Q}$, $k(W)= k(\d\cup S_1)$.  Hence,
$$n(W) = k(W) - v(G_{\Q}) + e(W) = k(\d\cup S_1) - k(\d)+ |S_1|.$$
\begin{eqnarray*}
n(\d\cup S_1)&=& k(\d\cup S_1)-v(\G)+ e(\d\cup S_1) \\
&=& \big(k(\d)-v(\G) + e(\d)\big) + \big(k(\d\cup S_1)-k(\d)+ |S_1|\big)\\
&=& n(\d) + n(W).
\end{eqnarray*}
%
Let $e\in\I(\Q)$. By Theorem \ref{treethm}, $\Q$ corresponds to $\r$
such that $e$ is nugatory, and by Lemma \ref{nugatorylemma},
$\Gamma_1(\r,e)$ is connected.
Since $f(\H)=|\gamma_{\H}|$, by Lemma \ref{nugatorylemma}, $f(\d\cup e)=f(\d)-1$.
Since live edges do not link each other, we can iterate this to obtain $f(\d\cup S_1)=f(\d)-|S_1|$.
Therefore, 
\begin{eqnarray*}
2 g(\d\cup S_1)&=& 2k(\d\cup S_1)-v(\G)+ e(\d\cup S_1) -f(\d\cup S_1) \\
&=&2k(\d)-v(\G)+\big(e(\d)+|S_1|\big)-\big(f(\d)-|S_1|\big) + 2k(\d\cup S_1)-2k(\d)\\
&=& 2g(\d)+2\big(k(\d\cup S_1)-k(\d)+|S_1|\big) \\
&=& 2g(\d) + 2 n(W).
\end{eqnarray*}
\done

{\em Proof (Theorem \ref{mainthm}): \quad } The sum in Theorem
\ref{mainthm} is over quasi-trees, which correspond to leaves $[\r]$ of
$\T$.  It suffices to show that for any quasi-tree, its summand in
Theorem \ref{mainthm} equals the sum over all $\H_s$ for $s\in[\r]$ in
equation (\ref{ssdbrt}).

Let $S = S_1 \cup S_2$, where $S_1 \subseteq \I$ and $S_2 \subseteq \E$.
By Lemma \ref{extlivecount}, the contribution from $[\r]$ to the sum in equation (\ref{ssdbrt}) is
\begin{eqnarray*}
& & \sum_{S \subseteq \I \cup \E } (X-1)^{k(\d\cup S)-1}\; Y^{n(\d\cup S)}\;Z^{g(\d\cup S)} \\
&=& \sum_{S_2 \subseteq \E} Y^{|S_2|} \sum_{S_1 \subseteq \I} (X-1)^{k(\d\cup S_1)-1}\; Y^{n(\d\cup S_1)}\;Z^{g(\d\cup S_1)} \\
&=& (1+Y)^{|\E|} \sum_{S_1 \subseteq \I} (X-1)^{k(\d\cup S_1)-1}\; Y^{n(\d\cup S_1)}\;Z^{g(\d\cup S_1)}.
\end{eqnarray*}

Below, we will use the spanning subgraph expansion of the
Tutte polynomial (see, e.g., p.339 of \c{ModernGraph}),
$$T_G(x,y)= \sum_{W\subseteq G} (x-1)^{k(W)-k(G)}\; (y-1)^{n(W)}.$$

Let $G_{\Q}$ denote the graph whose vertices are the components of
$\d$ and whose edges are the edges in $\I$.  $\Q$ is a connected
subgraph of $\G$, so $G_{\Q}$ is a connected graph, hence
$k(G_{\Q})=1$.  The subgraphs $\{\d\cup S_1\ |\ S_1 \subseteq \I\}$
are in one-one correspondence with spanning subgraphs $W\subseteq
G_{\Q}$.
Let $n_0= n(\d)$ and $g_0= g(\d)$.  By Lemma \ref{intlivecount},
\begin{eqnarray*}
&& \sum_{S_1 \subseteq \I} (X-1)^{k(\d\cup S_1)-1}\; Y^{n(\d\cup S_1)}\;Z^{g(\d\cup S_1)} \\
&=& \sum_{W \subseteq G_{\Q}} (X-1)^{k(W)-1}\; Y^{n(\d)+n(W)}\;Z^{g(\d)+n(W)} \\
&=& Y^{n_0}Z^{g_0}\sum_{W \subseteq G_{\Q}} (X-1)^{k(W)-k(G_{\Q})}\; (YZ)^{n(W)} \\
&=& Y^{n_0}Z^{g_0}\ T_{G_{\Q}}(X,1+YZ).
\end{eqnarray*}
The last step is obtained from the spanning subgraph expansion of the
Tutte polynomial with $x=X$ and $y=1+YZ$.  

Therefore, for each $\Q$, the contribution to the sum in (\ref{ssdbrt})
is 
$$ Y^{n(\d(\Q))}\;Z^{g(\d(\Q))}\;(1+Y)^{|\E(\Q)|}\;T_{G_{\Q}}(X,1+YZ). $$
This completes the proof of Theorem \ref{mainthm}.  
\done

\section{Example}
\label{sec_example}
We compute the quasi-tree and spanning tree expansions for a ribbon graph $\G$ with $12$
quasi-trees having a variety of topological types.  $\G$ has three vertices and six
edges, given by 
$\sigma_0  =  (1,3,2,5)\, (7, 9)\, (10,4,12,8,6,11)$, 
$\sigma_1 =  (1,2)\,(3,4)\,(5,6)\,(7,8)\,(9,10)\,(11,12)$, so 
$ \sigma_2  =  (1,6,7,10,12,3,2,4,9,8,11,5)$.
%
%
We order the edges of $\G$ by $\min(i,\sigma_1(i))$. The ribbon graph
$\G$ and its surface are shown below:

\begin{center}
\psfrag{1}{\footnotesize{$1$}}
\psfrag{2}{\footnotesize{$2$}}
\psfrag{3}{\footnotesize{$3$}}
\psfrag{4}{\footnotesize{$4$}}
\psfrag{5}{\footnotesize{$5$}}
\psfrag{6}{\footnotesize{$6$}}
\psfrag{7}{\footnotesize{$7$}}
\psfrag{8}{\footnotesize{$8$}}
\psfrag{9}{\footnotesize{$9$}}
\psfrag{10}{\footnotesize{$10$}}
\psfrag{11}{\footnotesize{$11$}}
\psfrag{12}{\footnotesize{$12$}}
\includegraphics[height=1in]{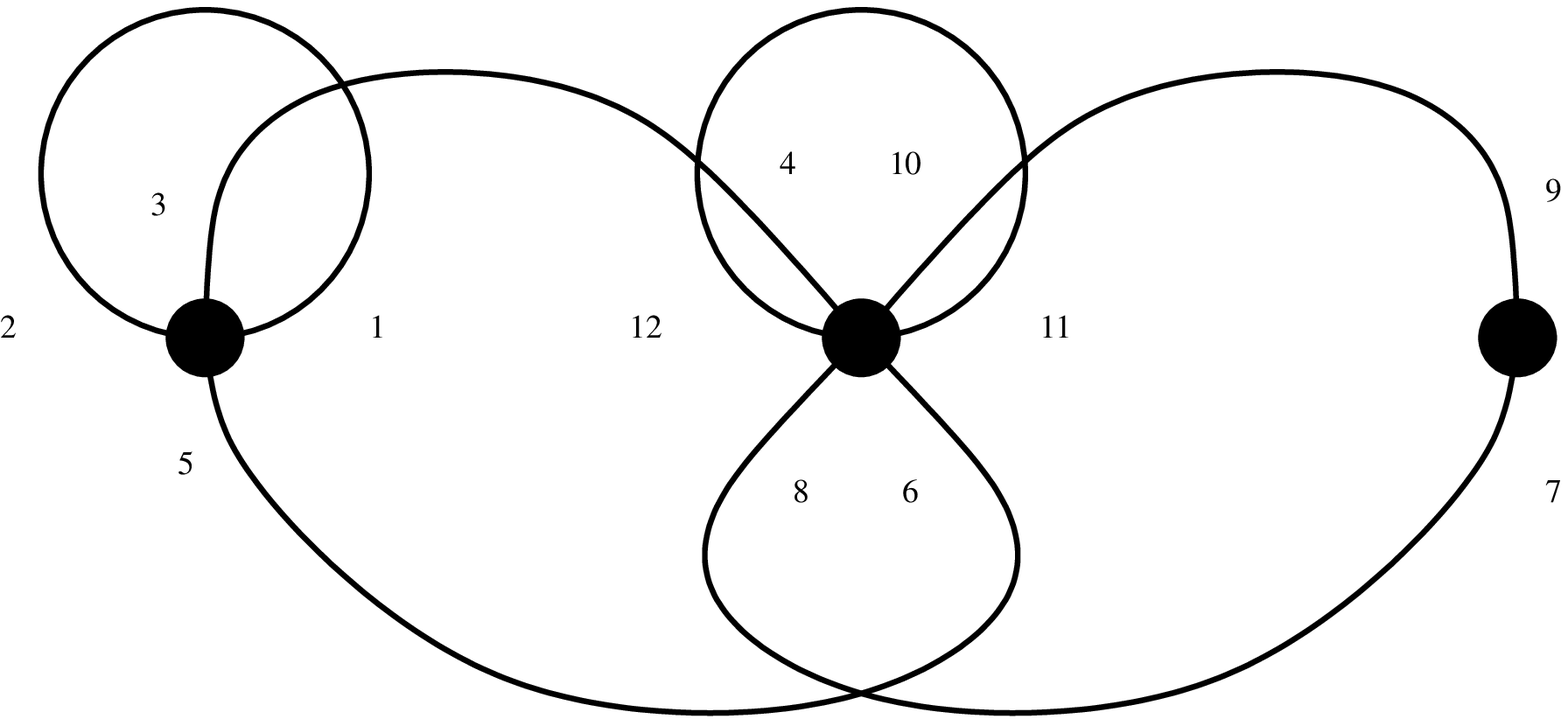}
\qquad
\includegraphics[height=0.9in]{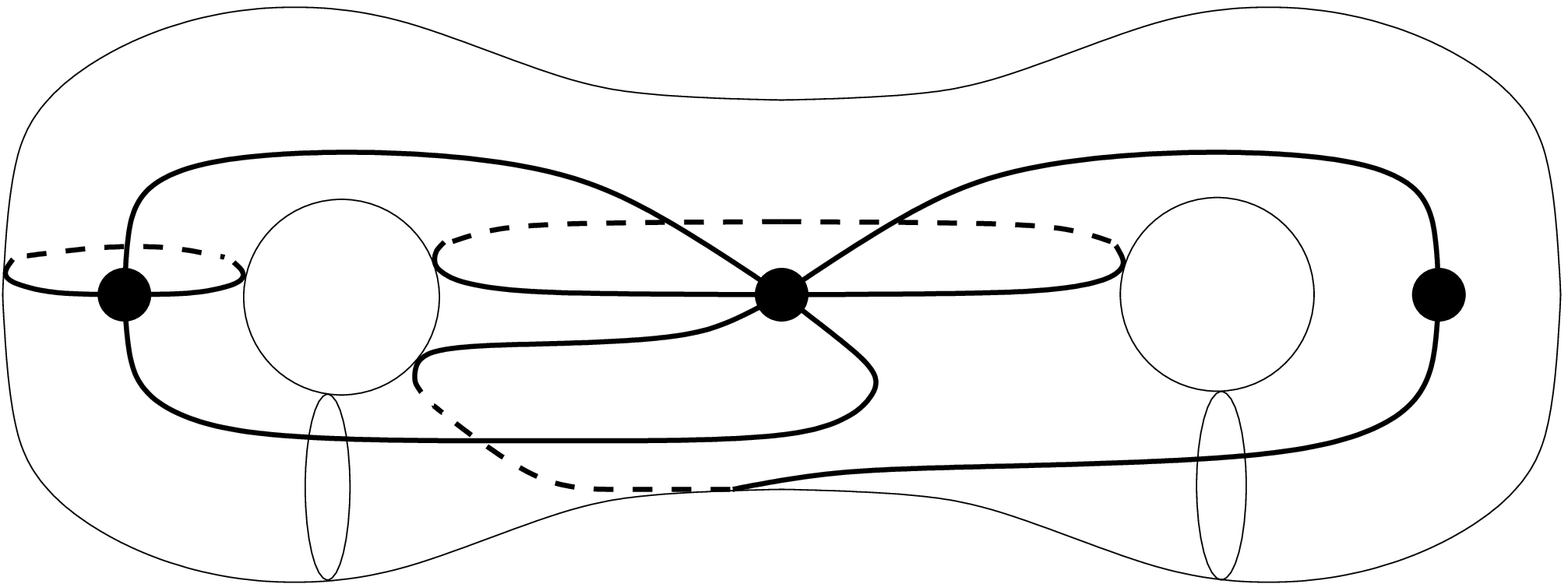}
 \end{center}

 In the table below, we denote quasi-trees using the edge order;
 e.g., 001010 denotes $\Q$ consisting of only the third and fifth
 edges, $(5,6)$ and $(9,10)$.  For each $\Q$, we compute the chord
 diagram, activites ($L,\ell$ for internally and externally live;
 $D,d$ for internally and externally dead), numbers $\{g,n,\bar g,\e\}
 = \{g(\Q), n(\d(\Q)), g(\d(\Q)), |\E(\Q)|\}$, graph $G_{\Q}$, and its
 weight in the sum of Theorem \ref{mainthm}.  For the chord diagrams,
 we give the cyclic permutation of the half-edges.  The types of
 graphs $G_{\Q}$ that occur in this example are as follows:

\begin{tabular}{ll}
1.  vertex & 2.  edge\\
3.  two edges with a vertex in common & 4.  two edges with both vertices in common\\
5.  2-cycle joined to a bridge & 6.  loop\\
7.  loop joined to a bridge. &
\end{tabular} \\
\scalebox{0.85}{\parbox{\textwidth}{%
\begin{center}
\begin{tabular}{|c|c|l|l|c|l|}
\hline
$\Q$ & $C_{\Q}$ & Activity & $g,n,\bar g,\e$ & $G_{\Q}$ & Weight \\
\hline
001010& $(1,3,2,5,11,10,7,9,4,12,8,5)$ & $\ell dDdDd$ & $0,0,0,1$ & $1$  & $(1+Y)$ \\
\hline
001100& $(1,3,2,5,11,10,4,12,8,9,7,6)$ & $\ell dDLdd$ & $0,0,0,1$ & $2$ & $X(1+Y)$ \\
\hline
001111& $(1,3,2,5,11,8,9,4,12,10,7,6)$& $\ell d DDDD$ & $1,2,1,1$ &$1$ & $Y^2Z(1+Y)$   \\
\hline
010010& $(1,3,12,8,6,11,10,7,9,4,2,5)$ & $\ell LddDd$ & $0,0,0,1$ & $2$ & $X(1+Y)$  \\
\hline
010100& $(1,3,12,8,9,7,6,11,10,4,2,5)$ & $\ell LdLdd$ & $0,0,0,1$ & $3$ & $X^2(1+Y)$ \\
\hline
010111& $(1,3,12,10,7,6,11,8,9,4,2,5)$ & $\ell LdDDD$ & $1,2,1,1$ & $2$ & $XY^2Z(1+Y)$ \\
\hline
011011& $(1,3,12,10,7,9,4,2,5,11,8,6)$ & $\ell LLdDD$ & $1,1,0,1$ & $4$ & $Y(1+Y)(X+1+YZ)$ \\
\hline
011101& $(1,3,12,10,4,2,5,11,8,9,7,6)$ & $\ell LLLdD$ & $1,1,0,1$ & $5$ & $XY(1+Y)(X+1+YZ)$ \\
\hline
011110& $(1,3,12,8,9,4,2,5,11,10,7,6)$ & $\ell LLDDd$ & $1,1,0,1$ & $4$ & $Y(1+Y)(X+1+YZ)$ \\
\hline
111010& $(1,5,11,10,7,9,4,2,3,12,8,6)$ &  $LDDdDd$ & $1,1,0,0$ & $6$ & $Y(1+YZ)$ \\
\hline
111100& $(1,5,11,10,4,2,3,12,8,9,7,6)$ & $LDDLdd$ & $1,1,0,0$ & $7$ & $XY(1+YZ)$ \\
\hline
111111& $(1,5,11,8,9,4,2,3,12,10,7,6)$ & $LDDDDD$ & $2,3,1,0$ & $6$ & $Y^3Z(1+YZ)$ \\
\hline
\end{tabular}
\end{center}
}}

Adding the weights in the last column, the \BRT polynomial of $\G$ is
\begin{eqnarray*}
C(\G)& = & Z^2 Y^4+2 X Z Y^3+4 Z Y^3+X^2 Y^2+3 X Y^2 + 3 X Z Y^2+4 Z Y^2+2 Y^2+
\\ & &
 2 X^2 Y+
 6 X Y+4 Y+X^2+2 X+1.
\end{eqnarray*}
By Proposition \ref{qtcount}, $q(\G;t,Y)=C(\G;1,Y,t Y^{-2})=4 + 7t + t^2$, which counts the quasi-trees of every genus.

As an example, let $\Q$ be the eighth quasi-tree, denoted 011101.  The
associated partial resolution is $\r=****01$. $\d(\Q)$ has three
components, consisting of two isolated vertices and a loop. $G_\Q$ has
three vertices and three edges, two connected in parallel and a second
edge to the remaining vertex.  The Tutte polynomial $T_{G_\Q}(x,y)=
x(x+y)$.  Thus, the contribution from $\Q$ is $XY(1+Y)(X+1+YZ)$, which
is also the contribution from the 16 terms in the state sum
(\ref{ssdbrt}) for all $s\in [\rho]$.

We now compute the spanning tree expansion (\ref{stbrt}) of the \BRT
polynomial for this example. Using the notation above, the spanning
trees are the genus zero quasi-trees. In the table below, we give the
spanning trees $T$, their activities in the sense of Tutte, their
weights given by the inner sum in (\ref{stbrt}), and the factor
$X^{i(T)}$ in (\ref{stbrt}).  Note that the activities for the
spanning trees below are different in every case from the activities
given above for the corresponding quasi-trees.

\begin{center}
\begin{tabular}{|c|l|c|l|}
\hline
$T$ & Activity & Weight & $X^{i(T)}$ \\
\hline
001010&  $\ell \ell D\ell D \ell$  & $1 + 4Y + 2Y^2 + 4Y^2Z + 4Y^3Z + Y^4Z^2 $ & 1\\
\hline
001100& $\ell \ell D Ld \ell$ & $1 + 3 Y + Y^2 + 2 Y^2
Z + Y^3 Z$ & $X$\\
\hline
010010& $\ell L d \ell D \ell$ & $(1 + Y)(1 + 2 Y + Y^2 Z)$ &$X$ \\
\hline
010100& $\ell LdLd\ell$ & $(1+Y)^2$ &$X^2$\\
\hline
\end{tabular}
\end{center}
Taking the sum according to (\ref{stbrt}), we obtain $C(\G)$ as above.

\bibliography{dkh}
\bibliographystyle{plain}

\end{document}